%% file: ci.tex
\documentclass[a4paper]{amsart}
\newcommand{\preamble}{
  \keywords{Fano manifolds, mirror symmetry, quantum differential equations, Picard--Fuchs equations}
  \email{t.coates@imperial.ac.uk}
  \email{a.m.kasprzyk@imperial.ac.uk}
  \email{t.prince12@imperial.ac.uk}
  \maketitle
}
\newcommand{\startintroblock}{}
\newcommand{\finishintroblock}{}
\newcommand{\noteondata}{The results of our computations in machine-readable form, together with full details of our implementation and all source code used, can be found in the ancillary files which accompany this paper on the {\tt arXiv}.}
\newcommand{\secondnoteondata}{The results of our computations, in machine readable form, together with full source code written in Magma~\cite{Magma} can be found in the ancillary files which accompany this paper on the {\tt arXiv}.  See the files called {\tt README.txt} for details. }
\newcommand{\Obrofootnote}{}
\usepackage[margin=2.45cm]{geometry}
\usepackage{hyperref}
\hypersetup{
    colorlinks,
    linkcolor={red!50!black},
    citecolor={blue!50!black},
    urlcolor={blue!80!black}
}
\usepackage{longtable}
\usepackage{pdflscape} % used by appendix
\usepackage{colortbl} % used by appendix
\usepackage{arydshln} % used by appendix
\usepackage{calc} % used by appendix
\newsavebox{\cimatrixbox}
\newlength{\cimatrixheight}
\newlength{\cimatrixshift}
\newcommand{\Obro}[2]{\text{\rm B\O S}^{#1}_{#2}} % The index in Obro's classification
\newcommand{\MM}[2]{\mathrm{MM}^3_{#1\text{--}#2}} % The Mori-Mukai number
\newcommand{\BB}[2]{B^{#1}_{#2}}
\newcommand{\VV}[2]{V^{#1}_{#2}}
\newcommand{\MW}[2]{\mathrm{MW}^{#1}_{#2}}
\renewcommand{\SS}[2]{S^2_{#2}}
\newcommand{\FI}[2]{\mathrm{FI}^{#1}_{#2}}
\newcommand{\cimatrix}[1]{%
  \tiny
    \begin{tabular}{c}%
    \savebox{\cimatrixbox}{#1}%
    \setlength{\cimatrixheight}{\totalheightof{\usebox{\cimatrixbox}}}%
    \setlength{\cimatrixshift}{0.5\cimatrixheight}%
    \addtolength{\cimatrixshift}{-0.5ex}%
    \raisebox{\cimatrixshift}{\usebox{\cimatrixbox}}%
    \addtolength{\cimatrixheight}{2.9ex}%
    \rule[-1.45ex]{0pt}{\cimatrixheight}%
    \end{tabular}%
}
\newcommand{\evnrow}{\rowcolor[gray]{0.95}}
\newcommand{\oddrow}{}
\usepackage{etex}
\usepackage{booktabs}
\usepackage[all]{xy}
\usepackage{cite}
\usepackage{tikz}
\usepackage{pgfplots}
\theoremstyle{plain}
\newtheorem{theorem}{\bf Theorem}[section]

\theoremstyle{definition}
\newtheorem{definition}[theorem]{\bf Definition}
\newtheorem{remark}[theorem]{\bf Remark}

\newtheorem{example}[theorem]{\bf Example}
\newcommand{\CC}{\mathbb{C}}

\newcommand{\Cstar}{\CC^\times}

\newcommand{\hG}{\widehat{G}}

\newcommand{\LL}{\mathbb{L}}
\DeclareMathOperator{\NC}{NC}

\newcommand{\cO}{\mathcal{O}}

\newcommand{\PP}{\mathbb{P}}
\DeclareMathOperator{\Pic}{Pic}
\newcommand{\QQ}{\mathbb{Q}}

\DeclareMathOperator{\rf}{rf}
\DeclareMathOperator{\rk}{rk}

\newcommand{\V}{\mathbb{V}}
\newcommand{\ZZ}{\mathbb{Z}}
\newcommand{\NumberOfTriples}{117173}
\newcommand{\NumberOfPeriods}{738}
\newcommand{\NumberOfNewPeriods}{527}
\begin{document}

\title{Four-Dimensional Fano Toric Complete Intersections}

\author{T. Coates, A. Kasprzyk and T. Prince}

\address{Department of Mathematics, Imperial College London, 180 Queen's Gate, London SW7 2AZ, UK}

\preamble

\begin{abstract}
  We find at least~527 new four-dimensional Fano manifolds, each of which is a complete intersection in a smooth toric Fano manifold.
\end{abstract}

\startintroblock
%--------------------------------------------------------------------------------
\section{Introduction}  \label{sec:introduction}
%--------------------------------------------------------------------------------

Fano manifolds are the basic building blocks of algebraic geometry, both in the sense of the Minimal Model Program~\cite{Reid,Corti,Birkar--Cascini--Hacon--McKernan,Hacon--McKernan} and as the ultimate source of most explicit examples and constructions.  There are finitely many deformation families of Fano manifolds in each dimension~\cite{Kollar--Miyaoka--Mori}.  There is precisely one $1$-dimensional Fano manifold: the line; there are 10 deformation families of $2$-dimensional Fano manifolds: the del~Pezzo surfaces; and there are 105~deformation families of $3$-dimensional Fano manifolds~\cite{Iskovskih:1,Iskovskih:2,Iskovskih:anticanonical,Mori--Mukai:Manuscripta,Mori--Mukai:Tokyo,Mori--Mukai:Kinosaki,Mori--Mukai:Erratum,Mori--Mukai:Turin}.  Very little is known about the classification of Fano manifolds in higher dimensions.  

In this paper we begin to explore the geography of Fano manifolds in dimension~$4$.   Four-dimensional Fano manifolds of higher Fano index have been classified~\cite{Shokurov,Kobayashi--Ochiai,Kollar,Serpico,Fujita:polarized_1,Fujita:polarized_2,Fujita:polarized_3,Fujita:book,Iskovskikh:Fano_1,Iskovskikh:anticanonical,Iskovskikh--Prokhorov}---there are~35 in total---but the most interesting case, where the Fano variety has index~1, is wide open.  We use computer algebra to find many 4-dimensional Fano manifolds that arise as complete intersections in toric Fano manifolds in codimension at most~$4$.  We find at least $\NumberOfPeriods$~examples, 717~of which have Fano index~1 and $\NumberOfNewPeriods$~of which are new.

Suppose that $Y$ is a toric Fano manifold and that $L_1,\ldots,L_c$ are nef line bundles on $Y$ such that ${-K_Y} - \Lambda$ is ample, where $\Lambda = c_1(L_1) + \cdots + c_1(L_c)$.  Let $X \subset Y$ be a smooth complete intersection defined by a regular section of $\oplus_i L_i$.  The Adjunction Formula gives that
\[
{-K_X} = \big({-K_Y}-\Lambda\big)\big|_X
\]
so $X$ is Fano.  We find all four-dimensional Fano manifolds $X$ of this form such that the codimension $c$ is at most $4$.  

Our interest in this problem is motivated by a program to classify Fano manifolds in higher dimensions using mirror symmetry \cite{CCGGK}.  For each $4$-dimensional Fano \finishintroblock manifold $X$ as above, therefore, we compute the essential ingredients for this program: the quantum period and regularized quantum differential equation associated to $X$, and a Laurent polynomial $f$ that corresponds to $X$ under mirror symmetry; we also calculate basic geometric data about $X$, the ambient space~$Y$, and~$f$.  \noteondata

%--------------------------------------------------------------------------------
\section{Finding Four-Dimensional Fano Toric Complete Intersections} \label{sec:finding}
%--------------------------------------------------------------------------------

Our method is as follows.  Toric Fano manifolds $Y$ are classified up to dimension~8 by Batyrev, Watanabe--Watanabe, Sato, Kreuzer--Nill, and \O{}bro\Obrofootnote.    For each toric Fano manifold~$Y$ of dimension $d=4+c$, we:
\begin{enumerate}
  \renewcommand{\theenumi}{\roman{enumi}}
\item compute the nef cone of $Y$;
\item find all $\Lambda \in H^2(Y;\ZZ)$ such that $\Lambda$ is nef and $-K_Y - \Lambda$ is ample;
\item decompose $\Lambda$ as the sum of $c$ nef line bundles $L_1,\ldots,L_c$ in all possible ways.
\end{enumerate}
Each such decomposition determines a 4-dimensional Fano manifold $X \subset Y$, defined as the zero locus of a regular section of the vector bundle $\oplus_i L_i$.  To compute the nef cone in step~(i), we consider dual exact sequences
\[
\xymatrix{
  0 \ar[rr] && 
  \LL \ar[rr] &&
  \ZZ^N \ar[rr]^{\rho} &&
  \ZZ^d \ar[rr]&& 
  0\\
  0  && 
  \LL^\star \ar[ll] &&
  (\ZZ^N)^\star \ar[ll]_D &&
  (\ZZ^d)^\star \ar[ll]_{\rho^\star} && 
  0 \ar[ll]}
\]
where the map $\rho$ is defined by the~$N$ rays of a fan $\Sigma$ for $Y$.  There are canonical identifications $\LL^\star \cong H^2(Y;\ZZ) \cong \Pic(Y)$, and the nef cone of $Y$ is the intersection of cones 
\[
\NC(Y) = \bigcap_{\sigma \in \Sigma} \left \langle D_i : i \not \in \sigma \right \rangle
\]
where $D_i$ is the image under $D$ of the $i$th standard basis vector in $(\ZZ^N)^\star$~\cite[Proposition~15.1.3]{Cox--Little--Schenck}.
%% to see this it suffices to show that the correspondence between chambers in H^2(Y;\RR) and fans Sigma is that a cone
%% sigma is in Sigma if and only if <D_i : i \not \in \sigma> covers the corresponding chamber.  But now Sigma is the normal
%% fan to the Delzant polytope for Y, so this is true. (see e.g. chapter VII of "Torus Actions on Symplectic Manifolds" by Audin)
The classes $\Lambda$ in step (ii) are the lattice points in the polyhedron $P = \NC(Y) \cap \big({-K_Y} - \NC(Y)\big)$ such that ${-K_Y} - \Lambda$ lies in the interior of $\NC(Y)$.  Since $\NC(Y)$ is a strictly convex cone, $P$ is compact and the number of lattice points in $P$ is finite.  We implement step~(iii) by first expressing $\Lambda$ as a sum of Hilbert basis elements in $\NC(Y)$ in all possible ways:
\begin{align}
  \label{eq:sum_of_basis_elements}
  \Lambda = b_1 + \cdots + b_r && \text{$b_i$ an element of the Hilbert basis for $\NC(Y)$}
\end{align}
where some of the $b_i$ may be repeated; this is a knapsack-style problem.  We then, for each decomposition \eqref{eq:sum_of_basis_elements}, partition the $b_i$ into $c$ subsets $S_1,\ldots,S_c$ in all possible ways, and define the line bundle $L_i$ to be the sum of the classes in $S_i$.  

We found~$\NumberOfTriples$ distinct triples $(X;Y;L_1,\ldots,L_c)$, with a total of 17934 distinct ambient toric varieties $Y$.  Note that the representation of a given Fano manifold $X$ as a toric complete intersection is far from unique:  for example, if $X$ is a complete intersection in $Y$ given by a section of $L_1\oplus\cdots\oplus L_c$ then it is also a complete intersection in $Y \times \PP^1$ given by a section of $\pi_1^\star L_1\oplus\cdots\oplus \pi_1^\star L_c \oplus \pi_2^\star \cO_{\PP^1}(1)$.  Thus we have found far fewer than~$\NumberOfTriples$ distinct four-dimensional Fano manifolds.  We show below, by calculating quantum periods of the Fano manifolds $X$, that we find at least~$\NumberOfPeriods$ non-isomorphic Fano manifolds.  Since the quantum period is a very strong invariant---indeed no examples of distinct Fano manifolds $X \not \cong X'$ with the same quantum period $G_X = G_{X'}$ are known---we believe that we found precisely~$\NumberOfPeriods$ non-isomorphic Fano manifolds.  Eliminating the quantum periods found in~\cite{known_4d}, we see that at least~$\NumberOfNewPeriods$ of our examples are new.

\begin{remark}
  There exist Fano manifolds which do not occur as complete intersections in toric Fano manifolds.  But in low dimensions, most Fano manifolds arise this way: 8~of the 10~del~Pezzo surfaces, and at least~78 of the~105 smooth 3-dimensional Fano manifolds, are complete intersections in toric Fano manifolds~\cite{QC105}.
\end{remark}
\begin{remark}
  It may be the case that any $d$-dimensional Fano manifold which occurs as a toric complete intersection in fact occurs as a toric complete intersection in codimension $d$; we know of no counterexamples.  But even if this holds in dimension~$4$, our search will probably not find all $4$-dimensional Fano manifolds which occur as toric complete intersections.  This is because, if one of the line bundles $L_i$ involved is nef but not ample, then the K\"ahler cone for $X$ can be strictly bigger than the K\"ahler cone for $Y$.  In other words, it is possible for ${-K_X}$ to be ample on $X$ even if ${-K_Y} - \Lambda$ is not ample on $Y$.  For an explicit example of this in dimension~$3$, see \cite[\S55]{QC105}.
\end{remark}

%--------------------------------------------------------------------------------
\section{Quantum Periods and Mirror Laurent Polynomials} \label{sec:quantum_periods}
%--------------------------------------------------------------------------------

The quantum period $G_X$  of a Fano manifold $X$ is a generating function
\begin{align}
  \label{eq:quantum_period}
  G_X(t) = 1 + \sum_{d=1}^\infty c_d t^d && \text{$t \in \CC$}
\end{align}
for certain genus-zero Gromov--Witten invariants~$c_d$ of $X$ which plays an important role in mirror symmetry.  A precise definition can be found in \cite[\S B]{QC105}, but roughly speaking one can think of $c_d$ as the `virtual number' of rational curves $C$ in $X$ that pass through a given point, satisfy certain constraints on their complex structure, and satisfy $\langle {-K_X}, C\rangle = d$.  The quantum period is discussed in detail in \cite{CCGGK,QC105}; for us what will be important is that the regularized quantum period
\begin{align}
  \label{eq:regularized_quantum_period}
  \hG_X(t) = 1 + \sum_{d=1}^\infty d! c_d t^d && \text{$t \in \CC$, $|t| \ll \infty$}
\end{align}
satisfies a differential equation called the \emph{regularized quantum differential equation} of $X$:
\begin{align}
  \label{eq:regularized_QDE}
  L_X \hG_X \equiv 0 && L_X = \sum_{m=0}^{m=N} p_m(t) D^m
\end{align}
where the $p_m$ are polynomials and $D = t \frac{d}{dt}$.

It has been proposed that Fano manifolds should correspond under mirror symmetry to Laurent polynomials which are \emph{extremal} or of \emph{low ramification}~\cite{CCGGK}, in the sense discussed in \S\ref{sec:ramification} below.  An $n$-dimensional Fano manifold $X$ is said to be \emph{mirror-dual} to a Laurent polynomial $f \in \CC[x_1^{\pm 1},\ldots,x_n^{\pm 1}]$ if the regularized quantum period of $X$ coincides with the classical period of $f$:
\begin{align*}
  \pi_f(t) & = \frac{1}{(2 \pi i)^n} 
  \int_{(S^1)^n} \frac{1}{1-tf} \frac{dx_1}{x_1} \cdots \frac{dx_n}{x_n} & \text{$t \in \CC$, $|t| \ll \infty$} 
\end{align*}
If a Fano manifold $X$ is mirror-dual to the Laurent polynomial $f$ then the regularized quantum differential equation \eqref{eq:regularized_QDE} for $X$ coincides with the Picard--Fuchs differential equation satisfied by $\pi_f$.  The correspondence between Fano manifolds and Laurent polynomials is not one-to-one, but it is expected that any two Laurent polynomials $f$,~$g$ that are mirror-dual to the same Fano manifold are related by a birational transformation $\varphi \colon (\Cstar)^n \dashrightarrow (\Cstar)^n$ called a \emph{mutation} or a \emph{symplectomorphism of cluster type}~\cite{mutations,Galkin--Usnich,Katzarkov--Przyjalkowski}: $\varphi^\star f = g$.  We will write such a mutation as $f \overset{\varphi}{\dashrightarrow} g$.  Mutations are known to preserve the classical period (ibid.): if $f \overset{\varphi}{\dashrightarrow} g$ then $\pi_f = \pi_g$.  

\begin{remark}
  \label{rem:the_PF}
  In the paragraphs above we discuss \emph{the} regularized quantum differential equation and \emph{the} Picard--Fuchs differential equation.  This involves choices of normalization. Our conventions are that the regularized quantum differential operator is the operator $L_X$ as in \eqref{eq:regularized_QDE} such that:
  \begin{enumerate}
    \renewcommand{\theenumi}{\roman{enumi}}
  \item \label{item:condition_1} the order, $N$, of $L_X$ is minimal; and
  \item the degree of $p_N(t)$ is minimal; and
  \item the leading coefficient of $p_N$ is positive; and
  \item \label{item:condition_4} the coefficients of the polynomials $p_0,\ldots,p_N$ are integers with greatest common divisor equal to $1$.
  \end{enumerate}
  The Picard--Fuchs differential operator is the differential operator $L_f$ such that:
  \begin{align*}
    L_f \,\pi_f \equiv 0 && L_f = \sum_{m=0}^{m=N} P_m(t) D^m
  \end{align*}
  where the $P_m$ are polynomials and $D = t \frac{d}{dt}$, and that the analogs of conditions \eqref{item:condition_1}--\eqref{item:condition_4} above hold.
\end{remark}

We determined the quantum period $G_X$, for each of the triples $(X;Y;L_1,\ldots,L_c)$ from~\S\ref{sec:finding}, as follows.  For each such triple we found, using the Mirror Theorem for toric complete intersections~\cite{Givental:toric} and a generalization of a technique due to V.~Przyjalkowski, a Laurent polynomial~$f$ that is mirror-dual to $X$.  This process is described in detail in \S\ref{sec:Przyjalkowski}.  We then computed, for each triple, the first 20 terms of the power series expansion of $\hG_X = \pi_f$ using the Taylor expansion:
\[
\pi_f(t) = \sum_{d=0}^\infty \alpha_d t^d 
\]
where $\alpha_d$ is the coefficient of the unit monomial in $f^d$. We divided the~$\NumberOfTriples$ triples into~$\NumberOfPeriods$ ``buckets'', according to the value of the first 20 terms of the power series expansion of $\hG_X = \pi_f$, and then proved that any two Fano manifolds $X$,~$X'$ in the same bucket have the same quantum period by exhibiting a chain of mutations $f \overset{\varphi_0}{\dashrightarrow} f_1 \overset{\varphi_1}{\dashrightarrow} \cdots \overset{\varphi_{n-1}}{\dashrightarrow} f_n \overset{\varphi_n}{\dashrightarrow} g$ that connects the Laurent polynomials $f$ and $g$ mirror-dual to $X$ and $X'$.  

For each quantum period $G_X$, we computed the quantum differential operator $L_X$ directly from the mirror Laurent polynomial $f$ chosen above, using Lairez's generalized Griffiths--Dwork algorithm~\cite{Lairez}.  The output from Lairez's algorithm is a differential operator $L = \sum_{m=0}^N P_m(t) D^m$ with $P_0,\ldots,P_N \in \QQ[t]$ such that, with very high probability, $L \pi_f \equiv 0$.  Such an operator $L$ gives a recurrence relation for the Taylor coefficients $\alpha_0$,~$\alpha_1$,~$\alpha_2$\ldots of $\pi_f$; using this recurrence relation and the first~$20$ Taylor coefficients computed above, we solved for the first~2000 Taylor coefficients $\alpha_k$.  We then consider an operator:
\[
\bar{L} = \sum_{m=0}^{\bar{N}} \bar{P}_m(t) D^m
\]
where the $\bar{P}_m$ are polynomials of degree at most $\bar{R}$, and impose the condition that $\bar{L} \pi_f \equiv 0$.  The~2000 Taylor coefficients of $\pi_f$ give 2000 linear equations for the coefficients of the polynomials $\bar{P}_m$ and, provided that $(\bar{N}+1)(\bar{R}+1) \ll  2000$, this linear system is highly over-determined.  Since we are looking for the Picard--Fuchs differential operator (see Remark~\ref{rem:the_PF}), we may assume that $(\bar{N},\bar{R})$ is lexicographically less than $(N,\deg p_N)$.  We searched systematically for such differential operators with $(\bar{N}+1)(\bar{R}+1) \ll 2000$, looking for the operator $\bar{L}$ with lexicographically minimal $(\bar{N},\bar{R})$ and clearing denominators so that the analogs of conditions (iii) and (iv) in Remark~\ref{rem:the_PF} holds.  We can say with high confidence that this operator $\bar{L}$ is in fact the Picard--Fuchs operator $L_f$, although this is not proven---partly because Lairez's algorithm relies on a randomized interpolation scheme that is not guaranteed to produce an operator annihilating $\pi_f$, and partly because if $L_f$ were to involve polynomials $P_m$ of extremely large degree, 2000~terms of the Taylor expansion of $\pi_f$ will not be enough to detect $L_f$.  The operators $\bar{L}$ that we found satisfy a number of delicate conditions that act as consistency checks: for example they are of Fuchsian type (which is true for $L_f$, as $L_f$ arises geometrically from a variation of Hodge structure).   Thus we are confident that $\bar{L} = L_f$ in every case\footnote{This could be proved in any given case using methods of van~Hoeij~\cite{van_Hoeij}; cf.~\cite[\S8.2.2]{Lairez}.}.  Since $\hG_X = \pi_f$ and $L_X = L_f$ by construction, this determines, with high confidence, the quantum period $G_X$ and the regularized quantum differential operator $L_X$.

\begin{remark}
  The use of Laurent polynomials and Lairez's algorithm is essential here.  There is a closed formula~\cite[Corollary~D.5]{QC105} for the quantum period of the Fano manifolds that we consider, and one could in principle use this together with the linear algebra calculation described above to compute (a good candidate for) the regularized quantum differential operator $L_X$.  In practice, however, for many of the examples that we treat here, it is impossible to determine enough Taylor coefficients from the formula: the  computations involved are well beyond the reach of current hardware, both in terms of memory consumption and runtime.  By contrast, our approach using mirror symmetry and Lairez's algorithm will run easily on a desktop PC.
\end{remark}

\begin{remark}
  The regularized quantum differential equation for $X$ coincides with the (unregularized) quantum differential equation for an anticanonical Calabi--Yau manifold $Z \subset X$.  The study of the regularized quantum period from this point of view was pioneered by Batyrev--Ciocan-Fontanine--Kim--van~Straten~\cite{BCFKvS:1,BCFKvS:2}, and an extensive study of fourth-order Calabi--Yau differential operators was made in~\cite{AESZ:paper}.  We found~26 quantum differential operators with $N=4$; these coincide with or are equivalent to the fourth-order Calabi--Yau differential operators with AESZ IDs 1, 3, 4, 5, 6, 15, 16, 17, 18, 19, 20, 21, 22, 23, 34, 369, 370, and 424 in the Calabi--Yau Operators Database~\cite{AESZ:database}, together with one new fourth-order Calabi--Yau differential operator (which corresponds to our period sequence with ID 469).
\end{remark}

%--------------------------------------------------------------------------------
\section{Ramification Data} \label{sec:ramification}
%--------------------------------------------------------------------------------

Consider now one of our regularized quantum differential operators:
\[
L_X = \sum_{m=0}^{m=N} p_m(t) D^m
\]
as in \eqref{eq:regularized_QDE}, and its local system $\V \to \PP^1 \setminus S$ of solutions.  Here $S \subset \PP^1$ is the set of singular points of the regularized quantum differential equation.

\begin{definition}[\!\!\protect{\cite{CCGGK}}]
  Let $S \subset \PP^1$ be a finite set and $\V \to \PP^1 \setminus S$ a local system.  Fix a basepoint $x \in \PP^1 \setminus S$.  For $s \in S$, choose a small loop that winds once anticlockwise around $s$ and connect it to $x$ via a path, thereby making a loop $\gamma_s$ about $s$ based at $x$.  Let $T_s \colon \V_x \to \V_x$ denote the monodromy of $\V$  along $\gamma_s$.  The \emph{ramification} of $\V$ is:
  \[
  \rf(\V) := \sum_{s \in S} \dim\Big(\V_x/{\V_x}^{\!\!\!T_s}\Big)
  \]
\end{definition}

The ramification $\rf(\V)$ is independent of the choices of basepoint $x$ and of small loops $\gamma_s$.  A non-trivial, irreducible local system $\V \to \PP^1 \setminus S$ has $\rf \V \ge 2 \rk \V$: see \cite[\S2]{CCGGK}.
\begin{definition}
  Let $\V \to \PP^1 \setminus S$ be a local system as above.  The \emph{ramification defect} of $\V$ is the quantity $\rf(\V) - 2\rk(\V)$.  A local system of ramification defect zero is called \emph{extremal}.  
\end{definition}

\begin{definition}
The \emph{ramification} (respectively \emph{ramification defect}) of a differential operator $L_X$ is the ramification (respectively ramification defect) of the local system of solutions $L_X f \equiv 0$.  
\end{definition}

To compute the ramification of $L_X$, we proceed as in~\cite{known_4d}.  One can compute  Jordan normal forms of the local $\log$-monodromies $\{\log T_s : s \in S\}$ using linear algebra over a splitting field $k$ for $p_N(t)$.  (Every singular point of $L_X$ is defined over $k$.)  This is classical, going back to Birkhoff~\cite{Birkhoff}, as corrected by Gantmacher~\cite[vol.~2,~\S10]{Gantmacher} and Turrittin~\cite{Turrittin}; a very convenient presentation can be found in the book of Kedlaya~\cite[\S7.3]{Kedlaya}.  In practice we use the symbolic implementation of $\overline{\QQ}$ provided by the computational algebra system Magma~\cite{Magma,Steel}.  We computed ramification data for~575 of the~738 regularized quantum differential operators, finding ramification defects as shown in Table~\ref{tab:ramification_defects}; this lends some support to the conjecture, due to Golyshev~\cite{CCGGK}, that a Laurent polynomial $f$ which is mirror-dual to a Fano manifold should have a Picard--Fuchs operator $L_f$ that is extremal or of low ramification.  For the remaining~163 regularized quantum differential operators, the symbol $p_N(t)$ contains a factor of extremely high degree.  This makes the computation of ramification data prohibitively expensive.

\begin{table}
  \begin{tabular}{rcccc} \toprule
    Ramification defect & 0 & 1 & 2 & 3 \\
    Number of occurrences & 92 & 290 & 167 & 26 \\
    \bottomrule \\
  \end{tabular}
  \caption{Ramification Defects for~575 of the~738 Regularized Quantum Differential Operators}
  \label{tab:ramification_defects}
\end{table}

%--------------------------------------------------------------------------------
\section{The Przyjalkowski Method} \label{sec:Przyjalkowski}
%--------------------------------------------------------------------------------

We now explain, given complete intersection data $(X;Y;L_1,\ldots,L_c)$ as in~\S\ref{sec:finding}, how to find a Laurent polynomial~$f$ that is mirror-dual to $X$.  This is a slight generalization of a technique that we learned from V.~Przyjalkowski\footnote{Przyjalkowski informs us that he learned this, for the case of the cubic threefold, from L.~Katzarkov~\cite{Katzarkov} and D.~Orlov~\cite{Przyjalkowski:personal_communication}.}~\cite{Przyjalkowski,Przyjalkowski:weak_LG}, and which is based on the mirror theorems for toric complete intersections due to Givental~\cite{Givental:toric} and Hori--Vafa~\cite{Hori--Vafa}.  Recall the exact sequence
\[
\xymatrix{
  0  && 
  \LL^\star \ar[ll] &&
  (\ZZ^N)^\star \ar[ll]_D &&
  (\ZZ^d)^\star \ar[ll]_{\rho^\star} && 
  0 \ar[ll]}
\]
from~\S\ref{sec:finding} and the elements $D_i \in \LL^\star$,~$1 \leq i \leq N$, defined by the standard basis elements of $(\ZZ^N)^\star$.  Recall further that $\LL^\star \cong \Pic(Y)$, so that each line bundle $L_m$ defines a class in  $\LL^\star$.  Suppose that there exists a choice of disjoint subsets $E$,~$S_1$,\ldots,$S_c$ of $\{1,2,\ldots,N\}$ such that:
\begin{itemize}
\item $\{D_j : j \in E\}$ is a basis for $\LL^\star$;
\item each $L_m$ is a non-negative linear combination of $\{D_j : j \in E\}$;
\item $\sum_{k \in S_m} D_k = L_m$ for each $m \in \{1,2,\ldots,c\}$;
\end{itemize}
and distinguished elements $s_m \in S_m$, $1 \leq m \leq c$.  Set $S_m^\circ = S_m \setminus \{s_m\}$.  Writing the map $D$ in terms of the standard basis for $(\ZZ^N)^\star$ and the basis $\{D_j : j \in E\}$ for $\LL^\star$ defines an  $(N-d) \times N$ matrix $(m_{ji})$ of integers.  Let $(x_1,\ldots,x_N)$ denote the standard co-ordinates on $(\Cstar)^N$, let $r = N-d$, and define $q_1,\ldots,q_r$ and $F_1,\ldots,F_c$ by:
\begin{align*}
  q_j = \prod_{i=1}^N x_i^{m_{ji}} &&
  F_m = \sum_{k \in S_m} x_k
\end{align*}
Givental~\cite{Givental:toric} and Hori--Vafa~\cite{Hori--Vafa} have shown that:
\begin{equation}
  \label{eq:mirror_theorem}
  G_X = \int_\Gamma e^{t W} \frac{\bigwedge_{i=1}^N \frac{dx_i}{x_i}}{\bigwedge_{m=1}^c d F_m \wedge \bigwedge_{j=1}^{r} \frac{dq_j}{q_j}}
\end{equation}
where $W = x_1 + \cdots + x_N$ and $\Gamma$ is a certain cycle in the submanifold of $(\Cstar)^N$ defined by:
\begin{align*}
  q_1 = \cdots = q_r = 1 && F_1 = \cdots = F_c = 1
\end{align*}
Introducing new variables $y_i$ for $i \in \bigcup_{m=1}^c S_m^\circ$, setting
\[
x_i = 
\begin{cases}
  \frac{1}{1 + \sum_{k \in S_m^\circ} y_k} & \text{if $i = s_m$} \\
  \frac{y_i}{1 + \sum_{k \in S_m^\circ} y_k} & \text{if $i \in S_m^\circ$} 
\end{cases}
\]
and using the relations $q_1 = \cdots = q_r = 1$ to eliminate the variables $x_j$,~$j \in E$, allows us to write $W - c$ as a Laurent polynomial $f$ in the variables:
\begin{align*}
  \left\{y_i : i \in \textstyle \bigcup_{m=1}^c S_m^\circ\right\}
  && \text{and} &&
  \left\{x_i : \text{$i \not \in E$ and $i \not \in \textstyle \bigcup_{m=1}^c S_m^\circ$}\right\}
\end{align*}
The mirror theorem~\eqref{eq:mirror_theorem} then implies that $\hG_X = \pi_f$, or in other words that $f$ is mirror-dual to $X$.

The Laurent polynomial $f$ produced by Przyjalkowski's method depends on our choices of $E$,~$S_1$,\ldots,$S_c$, and $s_1$,\ldots,$s_c$, but up to mutation this is not the case:

\begin{theorem}[\!\cite{Prince}]
  Let $Y$ be a toric Fano manifold and let $L_1,\ldots,L_c$ be nef line bundles on $Y$ such that ${-K_Y} - \Lambda$ is ample, where $\Lambda = c_1(L_1) + \cdots + c_1(L_c)$.  Let $X \subset Y$ be a smooth complete intersection defined by a regular section of $\oplus_i L_i$.  Let~$f$ and~$g$ be Laurent polynomial mirrors to $X$ obtained by applying Przyjalkowski's method to $(X;Y;L_1,\ldots,L_c)$ as above, but with possibly-different choices for the subsets $E$,~$S_1$,\ldots,$S_c$ and the elements $s_1$,\ldots,$s_c$.  Then there exists a mutation $\varphi$ such that $f \overset{\varphi}{\dashrightarrow} g$.
\end{theorem}

\begin{example}\label{eg:bundle}

Let $Y$ be the projectivization of the vector bundle $\cO^{\oplus 2} \oplus \cO(1)^{\oplus 2}$ over $\PP^2$.  Choose a basis for the two-dimensional lattice $\LL^\star$ such that the matrix $(m_{ji})$ of the map $D$ is:
\[
\begin{pmatrix}
1 & 1 & 1 & 0 & 0 & 1 & 1 \\
0 & 0 & 0 & 1 & 1 & 1 & 1
\end{pmatrix}
\]
Consider the line bundle $L_1 \to Y$ defined by the element $(2,1) \in \LL^\star$, and the Fano hypersurface $X \subset Y$ defined by a regular section of $L_1$.  Applying Przyjalkowski's method to the triple $(X;Y;L_1)$ with $E = \{3,4\}$, $S_1 = \{1,2,5\}$, and $s_1 = 1$ yields the Laurent polynomial
\[
f = \frac{(1+y_2+y_5)^2}{y_2 x_6 x_7} + \frac{1+y_2+y_5}{y_5 x_6 x_7} + x_6 + x_7
\]
mirror-dual to $X$.  Applying the method with $E = \{3,4\}$, $S_1 = \{1,6\}$, and $s_1 = 1$ yields:
\[
g  = x_2 + \frac{(1+y_6)^2}{x_2 y_6 x_7} + \frac{1+y_6}{x_5 y_6 x_7} + x_5 + x_7
\]
We have that $f \overset{\varphi}{\dashrightarrow} g$ where the mutation $\varphi \colon (\Cstar)^4 \to (\Cstar)^4$ is given by:
\[
(x_2,x_5,y_6,x_7) \mapsto \left( \frac{x_2}{x_5 y_6}, \frac{1}{y_6}, x_7, x_2+x_5 \right) = (y_2,y_5,x_6,x_7)
\]
\end{example}

\begin{remark}\label{rem:cod_dim}
Observe that, for a complete intersection of dimension~$n$ and codimension~$c$, Przyjalkowski's method requires partitioning $n+c$ variables into $c$ disjoint subsets.  
If $\frac{n+c}{c}<2$ then at least one of the subsets must have size one and so the corresponding variable, $x_j$ say, is eliminated from the Laurent polynomial via the equation $x_j = 1$.  One could therefore have obtained the resulting Laurent polynomial from a complete intersection with smaller codimension: new Laurent polynomials are found only when $\frac{n+c}{c} \geq 2$, that is, when the codimension is at most the dimension.  In particular, all possible mirrors to 4-dimensional Fano toric complete intersections given by the Przyjalkowski method occur for complete intersections in toric manifolds of dimension at most~8.
\end{remark}

%--------------------------------------------------------------------------------
\section{Examples}
%--------------------------------------------------------------------------------

\subsection{The Cubic 4-fold}

Let $X$ be the cubic $4$-fold.  This arises in our classification from the complete intersection data $(X; Y; L)$ with $Y = \PP^5$ and $L = \cO_{\PP^5}(3)$.  The Przyjalkowski method yields~\cite[\S2.1]{Ilten--Lewis--Przyjalkowski} a Laurent polynomial:
\[
f = \frac{(1+x+y)^3}{xyzw} + z + w
\]
mirror-dual to $X$, and elementary calculation gives:
\[
\pi_f(t) = \sum_{d=0}^\infty \frac{(3d)! (3d)!}{(d!)^6} t^{3d}
\]
Thus $\hG_X = \pi_f$, and the corresponding regularized quantum differential operator is:
\[
L_X = D^4 - 729 t^3 (D+1)^2(D+2)^2
\]
The local log-monodromies for the local system of solutions $L_X g \equiv 0$ are:
\begin{align*}
& \scriptsize \begin{pmatrix}
0&1&0&0\\
0&0&1&0\\
0&0&0&1\\
0&0&0&0
\end{pmatrix} && \text{at $\textstyle t=0$} \\& \scriptsize \begin{pmatrix}
0&0&0&0\\
0&0&0&0\\
0&0&0&1\\
0&0&0&0
\end{pmatrix} && \text{at $\textstyle t=\frac{1}{9}$} \\& \scriptsize \begin{pmatrix}
0&0&0&0\\
0&0&0&0\\
0&0&0&1\\
0&0&0&0
\end{pmatrix} && \text{at the roots of $\textstyle 81t^2+9t+1=0$} \\& \scriptsize \begin{pmatrix}
0&1&0&0\\
0&0&0&0\\
0&0&0&1\\
0&0&0&0
\end{pmatrix} && \text{at $t=\infty$} \\\end{align*}
and the operator $L_X$ is extremal.

\subsection{A (3,3) Complete Intersection in $\PP^6$}

Let $X$ be a complete intersection in $Y = \PP^6$ of type $(3,3)$.  This arises in our classification from the complete intersection data $(X; Y; L_1, L_2)$ with $L_1 = L_2 = \cO_{\PP^6}(3)$.  The Przyjalkowski method yields a Laurent polynomial:
\[
f = \frac{(1+x+y)^3(1+z+w)^3}{xyzw} -36
\]
mirror-dual to $X$, and~\cite[Corollary~D.5]{QC105} gives:
\[
\hG_X = \pi_f(t) = \sum_{k=0}^\infty \sum_{l=0}^\infty \frac{(3l)! (3l)! (k+l)!}{k! (l!)^7} (-36)^k t^{k+l}
\]
The corresponding regularized quantum differential operator $L_X$ is:
\begin{align*}
  (36t+1)^4(693t-1) & D^{4} \\
  + 18t(36t+1)^3(13860t+61) & D^{3}  \\
  + 9t(36t+1)^2(3492720t^2+57672t+77) & D^{2}  \\
  + 144t(36t+1)(11226600t^3+377622t^2+2754t+1) & D  \\
  + 15552t^2(1796256t^3+98496t^2+1605t+7) &
\end{align*}
The local log-monodromies for the local system of solutions $L_X g \equiv 0$ are:
\begin{align*}
& \scriptsize \begin{pmatrix}
0&1&0&0\\
0&0&1&0\\
0&0&0&1\\
0&0&0&0
\end{pmatrix} && \text{at $\textstyle t=0$} \\& \scriptsize \begin{pmatrix}
0&0&0&0\\
0&0&0&0\\
0&0&0&1\\
0&0&0&0
\end{pmatrix} && \text{at $\textstyle t=\frac{1}{693}$} \\& \scriptsize \begin{pmatrix}
\frac{2}{3}&1&0&0\\
0&\frac{2}{3}&0&0\\
0&0&\frac{1}{3}&1\\
0&0&0&\frac{1}{3}
\end{pmatrix} && \text{at $\textstyle t=-\frac{1}{36}$} \\\end{align*}
and so the operator $L_X$ is extremal.

%--------------------------------------------------------------------------------
\section{Results and Analysis}
%--------------------------------------------------------------------------------

\begin{figure}[htbp]
  \centering
  \includegraphics[width=0.9\textwidth]{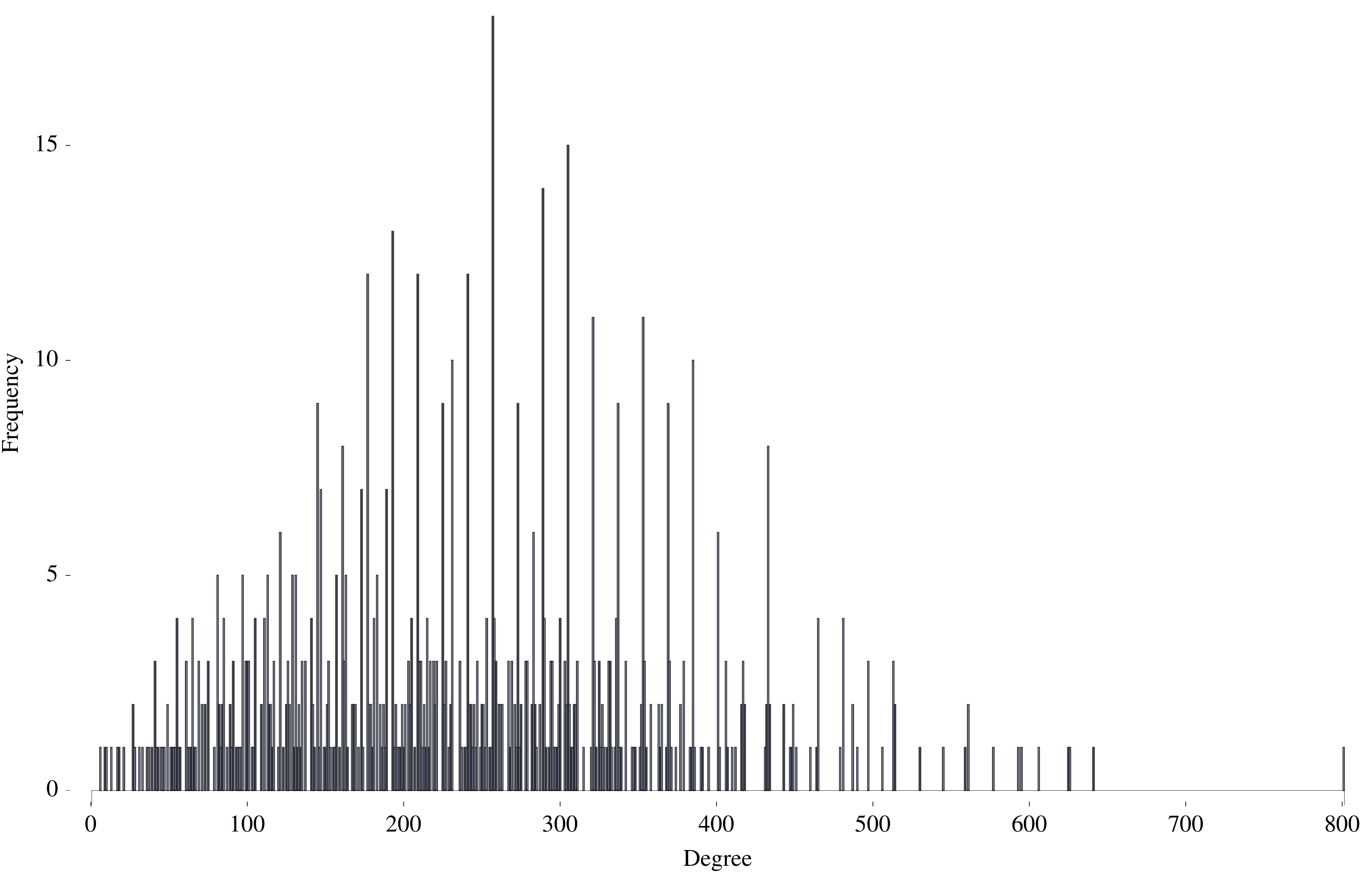}
  \caption{The Distribution of Degrees (Frequency Plot)}
  \label{fig:degree_bar_chart}
\end{figure}

\begin{figure}[htbp]
  \centering
  \includegraphics[width=0.9\textwidth]{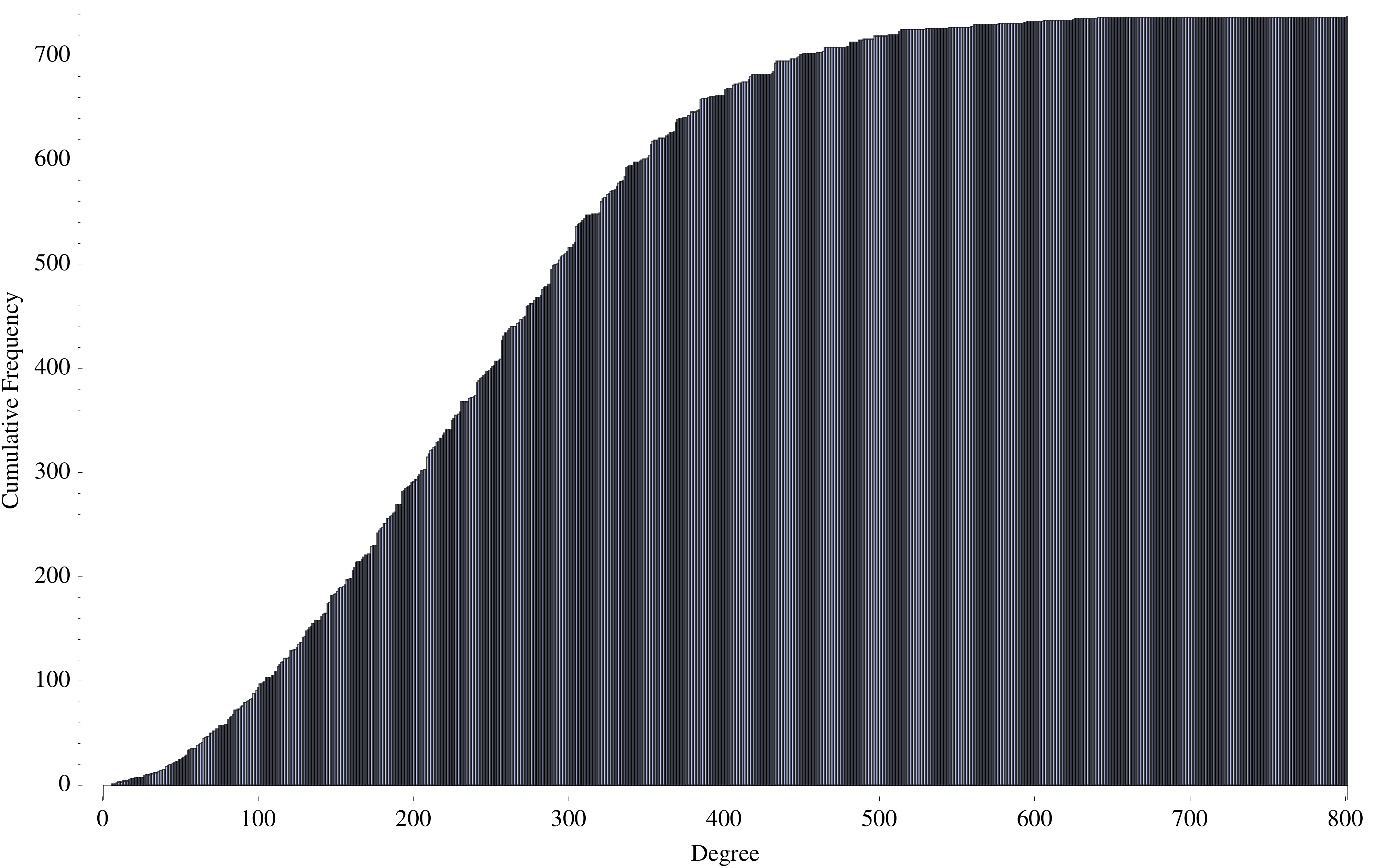}
  \caption{The Distribution of Degrees (Cumulative Frequency Plot)}
  \label{fig:degree_cumulative_frequency}
\end{figure}

We close by indicating how basic numerical invariants---degree and size of cohomology---vary across the 738 families of Fano manifolds that we have found.  The degree $({-K_X})^4$ varies from~$5$ to~$800$, as shown in Figures~\ref{fig:degree_bar_chart} and~\ref{fig:degree_cumulative_frequency}.

We do not have direct access to the size of the cohomology algebra of our Fano manifolds~$X$, as many of the line bundles occurring in the complete intersection data $(X;Y;L_1,\ldots,L_c)$ are not ample and so the Lefschetz Theorem need not apply.  But the order $N$ of the regularized quantum differential operator is a good proxy for the size of the cohomology. $N$ is the rank of a certain local system---an irreducible piece of the Fourier--Laplace transform of the restriction of the Dubrovin connection (in the Frobenius manifold given by the quantum cohomology of $X$) to the line in $H^\bullet(X)$ spanned by ${-K_X}$---and in the case where this local system is irreducible, which is typical, $N$ will coincide with the dimension of $H^\bullet(X)$.  For our examples, $N$ lies in the set $\{4,6,8,10,12\}$.  Figure~\ref{fig:degree_vs_N} shows how $N$ varies with the degree $({-K_X})^4$, with darker grays indicating a larger number of examples with that $N$ and degree. 

\pgfplotsset{
  compat=newest,
  every axis/.append style={font=\scriptsize,line width={0.1pt},tick style={major tick length=1pt,very thin,black}}
}
% we use buckets for degrees of size 10, with e.g. the bucket with degrees 1--10 and N=4 plotted as a square centred at (5,4), 
%                                                                             the bucket with degrees 11--20 and N=4 plotted as a square centred at (15,4), etc.
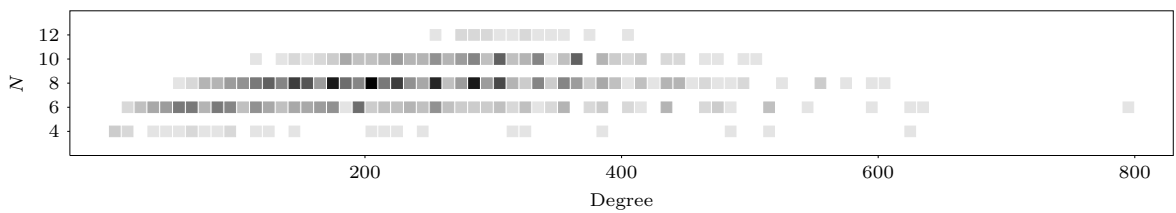
\begin{figure}[htbp]
  \centering
  \footnotesize
  \begin{tikzpicture}
    \begin{axis}[
      height=3.5cm,
      width=\textwidth,
      tick align=outside, unbounded coords=jump,
      ytick={4,6,8,10,12}, ytick pos=left,
      xtick={200,400,600,800}, xtick pos=left,
      xlabel=Degree, ylabel=$N$,
      xmin=-30, xmax=830, ymin=2, ymax=14,
      point meta min=0, point meta max=20, colormap={bw}{gray(0cm)=(0.94); gray(1cm)=(0)},
      scatter/use mapped color={draw=mapped color, fill = mapped color}]
      \addplot[mark=square*,only marks, scatter, scatter src=explicit,
      mark size=2]
      coordinates {
        (5,4) [3]
        (15,4) [2]
        (15,6) [2]
        (25,6) [4]
        (35,4) [1]
        (35,6) [6]
        (45,4) [1]
        (45,6) [7]
        (55,4) [1]
        (55,6) [10]
        (55,8) [1]
        (65,4) [2]
        (65,6) [10]
        (65,8) [2]
        (75,4) [1]
        (75,6) [5]
        (75,8) [5]
        (85,4) [1]
        (85,6) [10]
        (85,8) [5]
        (95,4) [2]
        (95,6) [9]
        (95,8) [7]
        (105,6) [4]
        (105,8) [8]
        (115,4) [1]
        (115,6) [8]
        (115,8) [10]
        (115,10) [1]
        (125,4) [1]
        (125,6) [6]
        (125,8) [12]
        (135,6) [4]
        (135,8) [9]
        (135,10) [1]
        (145,4) [1]
        (145,6) [6]
        (145,8) [15]
        (145,10) [2]
        (155,6) [6]
        (155,8) [13]
        (155,10) [1]
        (165,6) [7]
        (165,8) [7]
        (165,10) [2]
        (175,6) [8]
        (175,8) [18]
        (175,10) [3]
        (185,6) [1]
        (185,8) [11]
        (185,10) [6]
        (195,6) [11]
        (195,8) [9]
        (195,10) [4]
        (205,4) [1]
        (205,6) [3]
        (205,8) [20]
        (205,10) [4]
        (215,4) [1]
        (215,6) [4]
        (215,8) [10]
        (215,10) [5]
        (225,4) [1]
        (225,6) [4]
        (225,8) [15]
        (225,10) [7]
        (235,6) [5]
        (235,8) [8]
        (235,10) [5]
        (245,4) [1]
        (245,6) [4]
        (245,8) [6]
        (245,10) [5]
        (255,6) [8]
        (255,8) [17]
        (255,10) [8]
        (255,12) [1]
        (265,6) [4]
        (265,8) [4]
        (265,10) [5]
        (275,6) [3]
        (275,8) [7]
        (275,10) [7]
        (275,12) [2]
        (285,6) [3]
        (285,8) [18]
        (285,10) [9]
        (285,12) [2]
        (295,6) [3]
        (295,8) [7]
        (295,10) [4]
        (295,12) [2]
        (305,6) [4]
        (305,8) [14]
        (305,10) [12]
        (305,12) [1]
        (315,4) [1]
        (315,6) [2]
        (315,8) [5]
        (315,10) [4]
        (315,12) [1]
        (325,4) [1]
        (325,6) [4]
        (325,8) [3]
        (325,10) [5]
        (325,12) [2]
        (335,6) [1]
        (335,8) [9]
        (335,10) [9]
        (335,12) [1]
        (345,6) [2]
        (345,8) [3]
        (345,10) [1]
        (345,12) [1]
        (355,6) [5]
        (355,8) [9]
        (355,10) [4]
        (355,12) [1]
        (365,8) [7]
        (365,10) [12]
        (375,6) [2]
        (375,8) [3]
        (375,12) [1]
        (385,4) [1]
        (385,6) [3]
        (385,8) [6]
        (385,10) [5]
        (395,8) [4]
        (395,10) [3]
        (405,6) [2]
        (405,8) [1]
        (405,10) [2]
        (405,12) [1]
        (415,6) [1]
        (415,8) [4]
        (415,10) [3]
        (425,8) [1]
        (435,6) [5]
        (435,8) [5]
        (435,10) [2]
        (445,8) [5]
        (445,10) [2]
        (455,8) [1]
        (465,6) [2]
        (465,8) [2]
        (465,10) [1]
        (475,6) [3]
        (475,8) [1]
        (475,10) [1]
        (485,4) [1]
        (485,6) [1]
        (485,8) [1]
        (495,8) [2]
        (495,10) [1]
        (505,10) [1]
        (515,4) [1]
        (515,6) [4]
        (525,8) [1]
        (545,6) [1]
        (555,8) [3]
        (575,8) [1]
        (595,6) [1]
        (595,8) [1]
        (605,8) [1]
        (625,4) [1]
        (625,6) [1]
        (635,6) [1]
        (795,6) [1]
      };
    \end{axis}
  \end{tikzpicture}
  \caption{The Distribution of Degrees with $N$}
  \label{fig:degree_vs_N}
\end{figure}

\begin{figure}[htbp]
  \centering
  \footnotesize
  \begin{tikzpicture}
    \begin{axis}[
      height=3.5cm,
      width=\textwidth,
      tick align=outside, unbounded coords=jump,
      ytick={4,6,8,10,12}, ytick pos=left,
      xtick={200,400,600,800}, xtick pos=left,
      xlabel=Degree, ylabel=$N$,
      xmin=-30, xmax=830, ymin=2, ymax=14,
      point meta min=0, point meta max=20]
      \addplot[mark=square*,only marks, scatter, scatter src=explicit,
      colormap={bw}{gray(0cm)=(0.94); gray(1cm)=(0)},
      scatter/use mapped color={draw=mapped color, fill = mapped color},
      mark size=2]
      coordinates {
        (5,4) [3]
        (15,4) [2]
        (15,6) [2]
        (25,6) [4]
        (35,4) [1]
        (35,6) [6]
        (45,4) [1]
        (45,6) [7]
        (55,4) [1]
        (55,6) [10]
        (55,8) [1]
        (65,4) [2]
        (65,6) [10]
        (65,8) [2]
        (75,4) [1]
        (75,6) [5]
        (75,8) [5]
        (85,4) [1]
        (85,6) [10]
        (85,8) [5]
        (95,4) [2]
        (95,6) [9]
        (95,8) [7]
        (105,6) [4]
        (105,8) [8]
        (115,4) [1]
        (115,6) [8]
        (115,8) [10]
        (115,10) [1]
        (125,4) [1]
        (125,6) [6]
        (125,8) [12]
        (135,6) [4]
        (135,8) [9]
        (135,10) [1]
        (145,4) [1]
        (145,6) [6]
        (145,8) [15]
        (145,10) [2]
        (155,6) [6]
        (155,8) [13]
        (155,10) [1]
        (165,6) [7]
        (165,8) [7]
        (165,10) [2]
        (175,6) [8]
        (175,8) [18]
        (175,10) [3]
        (185,6) [1]
        (185,8) [11]
        (185,10) [6]
        (195,6) [11]
        (195,8) [9]
        (195,10) [4]
        (205,4) [1]
        (205,6) [3]
        (205,8) [20]
        (205,10) [4]
        (215,4) [1]
        (215,6) [4]
        (215,8) [10]
        (215,10) [5]
        (225,4) [1]
        (225,6) [4]
        (225,8) [15]
        (225,10) [7]
        (235,6) [5]
        (235,8) [8]
        (235,10) [5]
        (245,4) [1]
        (245,6) [4]
        (245,8) [6]
        (245,10) [5]
        (255,6) [8]
        (255,8) [17]
        (255,10) [8]
        (255,12) [1]
        (265,6) [4]
        (265,8) [4]
        (265,10) [5]
        (275,6) [3]
        (275,8) [7]
        (275,10) [7]
        (275,12) [2]
        (285,6) [3]
        (285,8) [18]
        (285,10) [9]
        (285,12) [2]
        (295,6) [3]
        (295,8) [7]
        (295,10) [4]
        (295,12) [2]
        (305,6) [4]
        (305,8) [14]
        (305,10) [12]
        (305,12) [1]
        (315,4) [1]
        (315,6) [2]
        (315,8) [5]
        (315,10) [4]
        (315,12) [1]
        (325,4) [1]
        (325,6) [4]
        (325,8) [3]
        (325,10) [5]
        (325,12) [2]
        (335,6) [1]
        (335,8) [9]
        (335,10) [9]
        (335,12) [1]
        (345,6) [2]
        (345,8) [3]
        (345,10) [1]
        (345,12) [1]
        (355,6) [5]
        (355,8) [9]
        (355,10) [4]
        (355,12) [1]
        (365,8) [7]
        (365,10) [12]
        (375,6) [2]
        (375,8) [3]
        (375,12) [1]
        (385,4) [1]
        (385,6) [3]
        (385,8) [6]
        (385,10) [5]
        (395,8) [4]
        (395,10) [3]
        (405,6) [2]
        (405,8) [1]
        (405,10) [2]
        (405,12) [1]
        (415,6) [1]
        (415,8) [4]
        (415,10) [3]
        (425,8) [1]
        (435,6) [5]
        (435,8) [5]
        (435,10) [2]
        (445,8) [5]
        (445,10) [2]
        (455,8) [1]
        (465,6) [2]
        (465,8) [2]
        (465,10) [1]
        (475,6) [3]
        (475,8) [1]
        (475,10) [1]
        (485,4) [1]
        (485,6) [1]
        (485,8) [1]
        (495,8) [2]
        (495,10) [1]
        (505,10) [1]
        (515,4) [1]
        (515,6) [4]
        (525,8) [1]
        (545,6) [1]
        (555,8) [3]
        (575,8) [1]
        (595,6) [1]
        (595,8) [1]
        (605,8) [1]
        (625,4) [1]
        (625,6) [1]
        (635,6) [1]
        (795,6) [1]
      };
      \addplot+[mark=square,only marks, scatter, scatter src=explicit,
      colormap={bw}{gray(0cm)=(0.94); gray(1cm)=(0)},
      scatter/use mapped color={draw=red, fill = mapped color},
      mark size=2, line width=0.5pt]
      coordinates {
        (215,4) [0]
        (225,4) [0]
        (245,6) [0]
        (255,8) [0]
        (265,8) [0]
        (275,10) [0]
        (285,6) [0]
        (285,8) [0]
        (285,12) [0]
        (295,6) [0]
        (295,8) [0]
        (295,12) [0]
        (305,6) [0]
        (305,8) [0]
        (305,10) [0]
        (305,12) [0]
        (315,12) [0]
        (325,6) [0]
        (325,8) [0]
        (325,10) [0]
        (325,12) [0]
        (335,8) [0]
        (335,10) [0]
        (335,12) [0]
        (345,10) [0]
        (345,12) [0]
        (355,8) [0]
        (355,10) [0]
        (355,12) [0]
        (365,8) [0]
        (365,10) [0]
        (375,8) [0]
        (375,12) [0]
        (385,4) [0]
        (385,6) [0]
        (385,8) [0]
        (385,10) [0]
        (395,8) [0]
        (395,10) [0]
        (405,6) [0]
        (405,8) [0]
        (405,10) [0]
        (405,12) [0]
        (415,8) [0]
        (415,10) [0]
        (435,6) [0]
        (435,8) [0]
        (435,10) [0]
        (445,8) [0]
        (445,10) [0]
        (455,8) [0]
        (465,6) [0]
        (465,8) [0]
        (465,10) [0]
        (475,6) [0]
        (475,8) [0]
        (475,10) [0]
        (485,4) [0]
        (485,6) [0]
        (485,8) [0]
        (495,8) [0]
        (495,10) [0]
        (505,10) [0]
        (515,6) [0]
        (525,8) [0]
        (545,6) [0]
        (555,8) [0]
        (575,8) [0]
        (595,6) [0]
        (595,8) [0]
        (605,8) [0]
        (625,4) [0]
        (635,6) [0]
        (795,6) [0]
      };
    \end{axis}
  \end{tikzpicture}
  \caption{The Distribution of Degrees with $N$, with Toric Fano Manifolds Highlighted.}
  \label{fig:degree_vs_N_with_toric}
\end{figure}
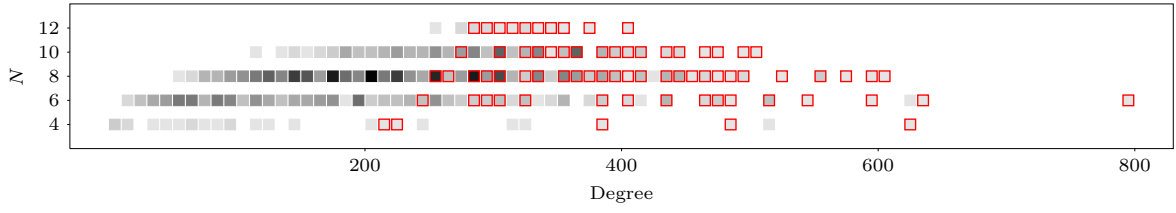

\begin{figure}[htbp]
  \centering
  \footnotesize
  %% bucketing in 10s, assigned to the middle, so anything in the range [10*k+1,10*k+10] gets plotted at 10*k+5
  \begin{tikzpicture}
    \begin{axis}[
      height=1.16279\textwidth,
      width=\textwidth,
      tick align=outside, unbounded coords=jump,
      ytick={-100,0,100,200,300,400,500,600,700,800},
      ytick pos=left,
      xtick={200,400,600,800}, xtick pos=left,
      xlabel=Degree, ylabel=Euler number,
      xmin=-30, xmax=830, ymin=-150, ymax=850,
      point meta min=0, point meta max=27]
      \addplot[mark=square*,only marks, scatter, scatter src=explicit,
      colormap={bw}{gray(0cm)=(0.94); gray(1cm)=(0)},
      scatter/use mapped color={draw=mapped color, fill = mapped color},
      mark size=2,line width=0pt]
      coordinates {
        (405,15) [6]
        (475,15) [4]
        (155,75) [1]
        (45,85) [1]
        (85,25) [1]
        (25,115) [1]
        (55,15) [1]
        (345,15) [6]
        (45,105) [1]
        (85,15) [1]
        (145,-5) [2]
        (345,25) [1]
        (85,55) [3]
        (145,25) [8]
        (305,25) [4]
        (45,115) [1]
        (515,5) [5]
        (55,105) [1]
        (85,35) [5]
        (85,45) [1]
        (55,85) [2]
        (145,15) [3]
        (305,15) [27]
        (55,95) [1]
        (145,5) [1]
        (45,25) [1]
        (55,65) [1]
        (55,75) [1]
        (145,45) [2]
        (145,35) [8]
        (115,65) [2]
        (595,15) [1]
        (45,-75) [1]
        (35,95) [1]
        (595,5) [1]
        (35,75) [1]
        (105,65) [1]
        (635,5) [1]
        (215,25) [6]
        (605,15) [1]
        (465,15) [5]
        (295,15) [12]
        (215,15) [13]
        (195,-15) [1]
        (185,35) [2]
        (5,825) [1]
        (225,15) [21]
        (125,65) [1]
        (355,15) [17]
        (155,-5) [1]
        (295,25) [4]
        (115,15) [1]
        (315,-5) [1]
        (225,25) [6]
        (175,-5) [1]
        (355,5) [2]
        (95,85) [1]
        (205,35) [2]
        (215,35) [1]
        (115,25) [6]
        (95,75) [1]
        (185,5) [1]
        (105,25) [3]
        (255,25) [5]
        (75,65) [2]
        (435,15) [8]
        (385,15) [15]
        (95,65) [1]
        (205,15) [17]
        (75,75) [2]
        (185,15) [6]
        (255,15) [26]
        (115,45) [2]
        (95,125) [1]
        (125,25) [7]
        (205,5) [2]
        (115,35) [5]
        (435,5) [4]
        (75,-25) [1]
        (165,115) [1]
        (185,25) [9]
        (255,5) [2]
        (5,365) [1]
        (105,15) [2]
        (415,15) [7]
        (205,25) [7]
        (65,45) [2]
        (105,55) [1]
        (65,35) [1]
        (165,65) [1]
        (125,5) [1]
        (75,35) [3]
        (115,55) [3]
        (415,5) [1]
        (235,5) [1]
        (445,15) [7]
        (15,215) [1]
        (95,-5) [1]
        (105,35) [5]
        (235,15) [11]
        (55,165) [1]
        (75,-5) [1]
        (75,45) [1]
        (95,25) [3]
        (95,-15) [1]
        (125,55) [1]
        (65,55) [3]
        (55,175) [1]
        (65,75) [1]
        (125,45) [3]
        (235,25) [4]
        (65,65) [2]
        (15,225) [1]
        (425,15) [1]
        (235,35) [1]
        (165,35) [2]
        (235,-15) [1]
        (5,555) [1]
        (95,-25) [1]
        (125,35) [4]
        (125,-5) [2]
        (395,15) [7]
        (55,135) [1]
        (65,95) [1]
        (35,195) [1]
        (95,55) [4]
        (65,85) [1]
        (525,15) [1]
        (65,-45) [1]
        (255,-5) [1]
        (625,5) [2]
        (175,45) [2]
        (75,25) [1]
        (325,15) [14]
        (25,255) [1]
        (175,35) [4]
        (165,15) [3]
        (95,35) [4]
        (575,15) [1]
        (155,35) [1]
        (335,15) [19]
        (545,5) [1]
        (325,25) [1]
        (175,55) [1]
        (505,15) [1]
        (35,155) [2]
        (175,15) [10]
        (165,25) [9]
        (15,295) [1]
        (335,25) [1]
        (35,145) [1]
        (155,25) [14]
        (175,5) [2]
        (195,25) [5]
        (115,-5) [1]
        (795,5) [1]
        (85,85) [1]
        (485,5) [1]
        (315,5) [2]
        (155,5) [1]
        (175,25) [9]
        (65,145) [1]
        (375,5) [1]
        (315,15) [10]
        (485,15) [2]
        (275,15) [16]
        (155,15) [2]
        (195,15) [14]
        (375,15) [5]
        (495,15) [3]
        (25,185) [1]
        (85,65) [1]
        (15,325) [1]
        (195,5) [1]
        (275,25) [3]
        (85,75) [2]
        (555,15) [3]
        (65,185) [1]
        (245,25) [6]
        (455,15) [1]
        (135,15) [2]
        (55,35) [1]
        (285,25) [5]
        (265,25) [3]
        (55,45) [1]
        (245,15) [10]
        (195,35) [3]
        (265,5) [1]
        (285,15) [26]
        (25,95) [1]
        (85,105) [1]
        (135,25) [7]
        (45,65) [1]
        (285,5) [1]
        (475,5) [1]
        (265,15) [9]
        (365,15) [19]
        (135,35) [5]
        (45,95) [2]
        (35,-115) [1]
      };
    \end{axis}
  \end{tikzpicture}
  \caption{The Distribution of Degrees with Euler Number.}
  \label{fig:degree_vs_Euler}
\end{figure}
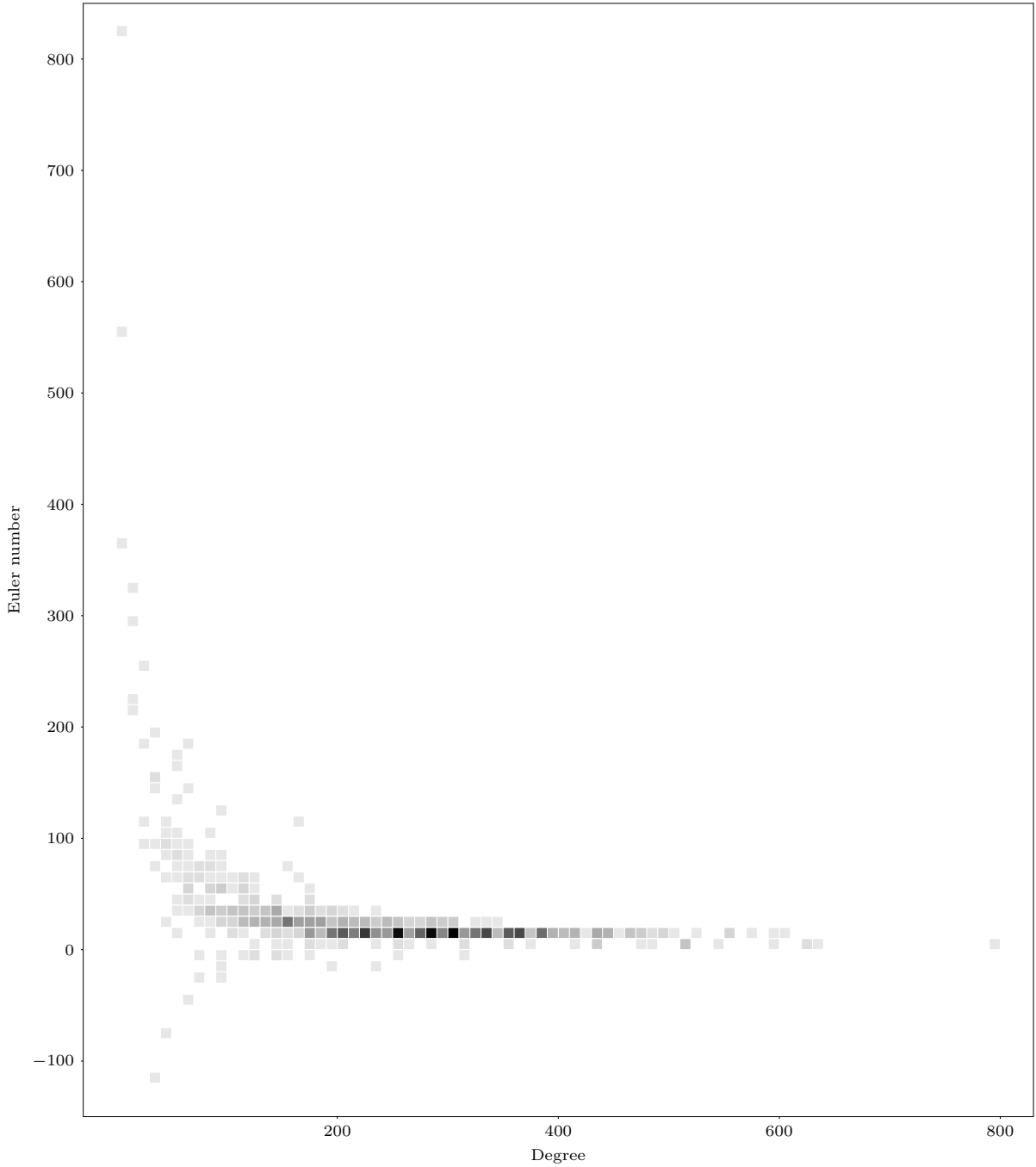

\noindent The isolated example on the right of Figure~\ref{fig:degree_vs_N}, with $N=6$ and degree~800, is the blow-up of $\PP(1,1,1,1,3)$ at a point.  Figure~\ref{fig:degree_vs_N_with_toric} again shows how $N$ varies with the degree $({-K_X})^4$, but this time with toric Fano manifolds highlighted in red.  Figure~\ref{fig:degree_vs_Euler} shows how the Euler number $\chi(T_X)$ varies with the degree $({-K_X})^4$, with darker grays indicating a larger number of examples with that Euler number and degree.  The three examples with the largest Euler number $\chi$ are a quintic hypersurface in $\PP^5$, with $\chi = 825$; a complete intersection of type $(2,4)$ in $\PP^6$, with $\chi = 552$; and a complete intersection of type $(3,3)$ in $\PP^6$, with $\chi = 369$.  The three examples with the most negative Euler number are $\PP^1 \times V^3_4$ where $V^3_4$ is a quartic hypersurface in $\PP^4$, with $\chi = {-112}$; $\PP^1 \times V^3_6$ where $V^3_6$ is a complete intersection of type $(2,3)$ in $\PP^5$, with $\chi = {-72}$; and $\PP^1 \times V^3_8$ where $V^3_8$ is a complete intersection of type $(2,2,2)$ in $\PP^6$, with $\chi = {-48}$.

%--------------------------------------------------------------------------------
\section{Source Code and Data}
%--------------------------------------------------------------------------------

\secondnoteondata The source code and data, but not the text of this paper, are released under a Creative Commons~CC0 license~\cite{CC0}: see the files called {\tt COPYING.txt} for details.  If you make use of the source code or data in an academic or commercial context, you should acknowledge this by including a reference or citation to this paper.  

%--------------------------------------------------------------------------------
\section*{Acknowledgments}
%--------------------------------------------------------------------------------

The computations underlying this work were performed using the Imperial College High Performance Computing Service and the compute cluster at the Department of Mathematics, Imperial College London.  We thank Simon Burbidge, Matt Harvey, and Andy Thomas for valuable technical assistance.  We thank John Cannon and the Computational Algebra Group at the University of Sydney for providing licenses for the computer algebra system Magma.  This research was supported by a Royal Society University Research Fellowship (TC); the Leverhulme Trust; ERC Starting Investigator Grant number~240123; EPSRC grant EP/I008128/1; an EPSRC Small Equipment grant; and an EPSRC Prize Studentship (TP).  We thank Alessio Corti for a number of very useful conversations, Pierre Lairez for explaining his generalized Griffiths--Dwork algorithm and sharing his code with us, and Duco van Straten for his analysis of the regularized quantum differential operators with $N=4$.

%--------------------------------------------------------------------------------
\bibliographystyle{amsplain}
\bibliography{bibliography}
%--------------------------------------------------------------------------------

%--------------------------------------------------------------------------------
% appendix
%--------------------------------------------------------------------------------
\pagebreak
\appendix
\input{appendix.tex}
\end{document}

%% file: appendix.tex
\begin{landscape}
\section{Regularized Quantum Period Sequences Arising From Complete Intersections in Toric Fano Manifolds}
\label{appendix:periods}

In this Appendix we record, for each of the toric complete intersections $X$ considered in this paper, the description, degree, and first few terms of the regularized quantum period sequence for $X$, as well as a representative construction of $X$ as a toric complete intersection.  The regularized quantum period sequences, in lexicographic order, are shown in Table~\ref{tab:toric_ci_periods}.  If the description there is left blank then no Fano manifold with that regularized quantum period sequence was previously known; otherwise the descriptions are exactly as in~\cite[Appendix~A]{known_4d}.  The quantities $\alpha_0,\alpha_1,\ldots$ are Taylor coefficients of the regularized quantum period sequence: $\hG_X(t) = \sum_{d=0}^\infty \alpha_d t^d$.  Ten terms of the Taylor expansion suffices to distinguish all the regularized quantum periods in Table~\ref{tab:toric_ci_periods}: see the discussion and supplementary table on page~\hyperlink{tab:extra_terms}{\pageref*{pages:extra_terms}}.  Tables~\ref{tab:index_1},~\ref{tab:index_2},~\ref{tab:index_3},~\ref{tab:index_4}, and~\ref{tab:index_5} contain constructions of our Fano manifolds~$X$ as toric complete intersections, together with the description, degree, Euler number, and first two terms of the Hilbert series for $X$.

\setlength{\LTcapwidth}{24cm}
\setlength{\extrarowheight}{0.2em}
{ \small
\input{table_ps}
}
\pagebreak
\end{landscape}
%-------------------------------------------------------------------------------
\newgeometry{left=1cm,right=1cm,top=2.45cm,bottom=2.45cm}

% Extra terms needed for period ids [ 81, 82, 117, 118, 248, 249 ]
It appears from Table~\ref{tab:toric_ci_periods} as if the regularized quantum period might coincide for the three pairs with period IDs $81$ and $82$, $117$ and $118$, and $248$ and $249$. This is not the case. The coefficients $\alpha_8$,~$\alpha_9$, and~$\alpha_{10}$ in these cases are:\\
\begin{center}
\small
\label{pages:extra_terms}
\setlength{\extrarowheight}{0.2em}
\begin{tabular}{cccc}
\toprule
Period ID \hypertarget{tab:extra_terms} &$\alpha_{8}$&$\alpha_{9}$&$\alpha_{10}$\\
\midrule
\oddrow $81$&$10990$&$102480$&$365652$\\
\evnrow $82$&$14350$&$87360$&$441252$\\
\oddrow $117$&$32830$&$212520$&$1190952$\\
\evnrow $118$&$32830$&$227640$&$1190952$\\
\oddrow $248$&$120820$&$814800$&$5670504$\\
\evnrow $249$&$127540$&$814800$&$6048504$\\
\bottomrule
\end{tabular}
\end{center}
\vspace{2em}
%-------------------------------------------------------------------------------
\small
\input{table_1}
%-------------------------------------------------------------------------------
\input{table_2}
%-------------------------------------------------------------------------------
\input{table_3}
%-------------------------------------------------------------------------------
\input{table_4}
%-------------------------------------------------------------------------------
\setlength{\LTcapwidth}{\textwidth}
\input{table_5}
\vspace{3em}

%% file: table_ps.tex
\begin{longtable}{ccccccccccc}
\caption{$738$~regularized period sequences obtained from $4$-dimensional Fano manifolds that arise as complete intersections in toric Fano manifolds.}
\label{tab:toric_ci_periods}\\
\toprule
\multicolumn{1}{c}{Period ID}&\multicolumn{1}{c}{Description}&\multicolumn{1}{c}{$(-K_X)^4$}&\multicolumn{1}{c}{$\alpha_{0}$}&\multicolumn{1}{c}{$\alpha_{1}$}&\multicolumn{1}{c}{$\alpha_{2}$}&\multicolumn{1}{c}{$\alpha_{3}$}&\multicolumn{1}{c}{$\alpha_{4}$}&\multicolumn{1}{c}{$\alpha_{5}$}&\multicolumn{1}{c}{$\alpha_{6}$}&\multicolumn{1}{c}{$\alpha_{7}$}\\
\midrule
\endfirsthead
\multicolumn{11}{l}{\tiny Continued from previous page.}\\
\addlinespace[1.7ex]
\midrule
\multicolumn{1}{c}{Period ID}&\multicolumn{1}{c}{Description}&\multicolumn{1}{c}{$(-K_X)^4$}&\multicolumn{1}{c}{$\alpha_{0}$}&\multicolumn{1}{c}{$\alpha_{1}$}&\multicolumn{1}{c}{$\alpha_{2}$}&\multicolumn{1}{c}{$\alpha_{3}$}&\multicolumn{1}{c}{$\alpha_{4}$}&\multicolumn{1}{c}{$\alpha_{5}$}&\multicolumn{1}{c}{$\alpha_{6}$}&\multicolumn{1}{c}{$\alpha_{7}$}\\
\midrule
\endhead
\midrule
\multicolumn{11}{r}{\tiny Continued on next page.}\\
\endfoot
\bottomrule
\endlastfoot
\oddrow \hypertarget{ps:1}{}\hyperlink{desc:1}{$1$}&$\PP^4$&$625$&$1$&$0$&$0$&$0$&$0$&$120$&$0$&$0$\\
\evnrow \hypertarget{ps:2}{}\hyperlink{desc:2}{$2$}&$\Obro{4}{115}$&$512$&$1$&$0$&$0$&$0$&$24$&$120$&$0$&$0$\\
\oddrow \hypertarget{ps:3}{}\hyperlink{desc:3}{$3$}&$Q^4$&$512$&$1$&$0$&$0$&$0$&$48$&$0$&$0$&$0$\\
\evnrow \hypertarget{ps:4}{}\hyperlink{desc:4}{$4$}&&$431$&$1$&$0$&$0$&$0$&$48$&$120$&$0$&$0$\\
\oddrow \hypertarget{ps:5}{}\hyperlink{desc:5}{$5$}&&$376$&$1$&$0$&$0$&$0$&$72$&$120$&$0$&$0$\\
\evnrow \hypertarget{ps:6}{}\hyperlink{desc:6}{$6$}&&$341$&$1$&$0$&$0$&$0$&$96$&$120$&$0$&$0$\\
\oddrow \hypertarget{ps:7}{}\hyperlink{desc:7}{$7$}&$\Obro{4}{21}$&$594$&$1$&$0$&$0$&$6$&$0$&$0$&$90$&$1260$\\
\evnrow \hypertarget{ps:8}{}\hyperlink{desc:8}{$8$}&$\Obro{4}{118}$&$513$&$1$&$0$&$0$&$6$&$0$&$120$&$90$&$0$\\
\end{longtable}

%% file: table_1.tex
\begin{longtable}{c@{}c@{}c@{}c@{}c@{}c@{}c}
\caption{Certain $4$-dimensional Fano manifolds with Fano index $r=1$ that arise as complete intersections in toric Fano manifolds.}
\label{tab:index_1}\\
\toprule
\multicolumn{1}{c}{Period ID}&\multicolumn{1}{c}{Description}&\multicolumn{1}{c}{$(-K_X)^4$}&\multicolumn{1}{c}{$\chi(T_X)$}&\multicolumn{1}{c}{$h^0(-K_X)$}&\multicolumn{1}{c}{$h^0(-2K_X)$}&\multicolumn{1}{c}{Weights and line bundles}\\
\midrule
\endfirsthead
\multicolumn{7}{l}{\tiny Continued from previous page.}\\
\addlinespace[1.7ex]
\midrule
\multicolumn{1}{c}{Period ID}&\multicolumn{1}{c}{Description}&\multicolumn{1}{c}{$(-K_X)^4$}&\multicolumn{1}{c}{$\chi(T_X)$}&\multicolumn{1}{c}{$h^0(-K_X)$}&\multicolumn{1}{c}{$h^0(-2K_X)$}&\multicolumn{1}{c}{Weights and line bundles}\\
\midrule
\endhead
\midrule
\multicolumn{7}{r}{\tiny Continued on next page.}\\
\endfoot
\bottomrule
\endlastfoot
\oddrow \hypertarget{desc:738}{}\hyperlink{ps:738}{$738$}&&$5$&$825$&$6$&$21$&\cimatrix{$\begin{array}{cccccc|c}
1&1&1&1&1&1&5\\\end{array}$}\\
\evnrow \hypertarget{desc:737}{}\hyperlink{ps:737}{$737$}&&$8$&$552$&$7$&$27$&\cimatrix{$\begin{array}{ccccccc|cc}
1&1&1&1&1&1&1&2&4\\\end{array}$}\\
\oddrow \hypertarget{desc:736}{}\hyperlink{ps:736}{$736$}&&$9$&$369$&$7$&$28$&\cimatrix{$\begin{array}{ccccccc|cc}
1&1&1&1&1&1&1&3&3\\\end{array}$}\\
\evnrow \hypertarget{desc:735}{}\hyperlink{ps:735}{$735$}&&$12$&$324$&$8$&$34$&\cimatrix{$\begin{array}{cccccccc|ccc}
1&1&1&1&1&1&1&1&2&2&3\\\end{array}$}\\
\oddrow \hypertarget{desc:731}{}\hyperlink{ps:731}{$731$}&&$16$&$224$&$9$&$41$&\cimatrix{$\begin{array}{ccccccccc|cccc}
1&1&1&1&1&1&1&1&1&2&2&2&2\\\end{array}$}\\
\evnrow \hypertarget{desc:734}{}\hyperlink{ps:734}{$734$}&&$17$&$293$&$10$&$45$&\cimatrix{$\begin{array}{ccccccc|c}
1&1&0&0&0&0&0&1\\0&0&1&1&1&1&1&4\\\end{array}$}\\
\oddrow \hypertarget{desc:733}{}\hyperlink{ps:733}{$733$}&&$20$&$212$&$11$&$51$&\cimatrix{$\begin{array}{ccccccc|c}
1&1&0&1&1&1&1&4\\0&0&1&1&1&1&1&4\\\end{array}$}\\
\evnrow \hypertarget{desc:730}{}\hyperlink{ps:730}{$730$}&&$26$&$186$&$12$&$60$&\cimatrix{$\begin{array}{cccccccc|cc}
1&1&0&0&0&0&0&0&0&1\\0&0&1&1&1&1&1&1&2&3\\\end{array}$}\\
\end{longtable}

%% file: table_2.tex
\begin{longtable}{ccccccc}
\caption{Certain $4$-dimensional Fano manifolds with Fano index $r=2$ that arise as complete intersections in toric Fano manifolds.}
\label{tab:index_2}\\
\toprule
\multicolumn{1}{c}{Period ID}&\multicolumn{1}{c}{Description}&\multicolumn{1}{c}{$(-K_X)^4$}&\multicolumn{1}{c}{$\chi(T_X)$}&\multicolumn{1}{c}{$h^0(-K_X)$}&\multicolumn{1}{c}{$h^0(-2K_X)$}&\multicolumn{1}{c}{Weights and line bundles}\\
\midrule
\endfirsthead
\multicolumn{7}{l}{\tiny Continued from previous page.}\\
\addlinespace[1.7ex]
\midrule
\multicolumn{1}{c}{Period ID}&\multicolumn{1}{c}{Description}&\multicolumn{1}{c}{$(-K_X)^4$}&\multicolumn{1}{c}{$\chi(T_X)$}&\multicolumn{1}{c}{$h^0(-K_X)$}&\multicolumn{1}{c}{$h^0(-2K_X)$}&\multicolumn{1}{c}{Weights and line bundles}\\
\midrule
\endhead
\midrule
\multicolumn{7}{r}{\tiny Continued on next page.}\\
\endfoot
\bottomrule
\endlastfoot
\oddrow \hypertarget{desc:682}{}\hyperlink{ps:682}{$682$}&$\VV{4}{4}$&$64$&$188$&$21$&$125$&\cimatrix{$\begin{array}{cccccc|c}
1&1&1&1&1&1&4\\\end{array}$}\\
\evnrow \hypertarget{desc:640}{}\hyperlink{ps:640}{$640$}&$\VV{4}{6}$&$96$&$90$&$27$&$175$&\cimatrix{$\begin{array}{ccccccc|cc}
1&1&1&1&1&1&1&2&3\\\end{array}$}\\
\oddrow \hypertarget{desc:596}{}\hyperlink{ps:596}{$596$}&$\VV{4}{8}$&$128$&$48$&$33$&$225$&\cimatrix{$\begin{array}{cccccccc|ccc}
1&1&1&1&1&1&1&1&2&2&2\\\end{array}$}\\
\evnrow \hypertarget{desc:567}{}\hyperlink{ps:567}{$567$}&$\PP^1\times \BB{3}{3}$&$192$&$-12$&$45$&$325$&\cimatrix{$\begin{array}{ccccccc|c}
1&1&0&0&0&0&0&0\\0&0&1&1&1&1&1&3\\\end{array}$}\\
\oddrow \hypertarget{desc:500}{}\hyperlink{ps:500}{$500$}&$\PP^1\times \BB{3}{4}$&$256$&$0$&$57$&$425$&\cimatrix{$\begin{array}{cccccccc|cc}
1&1&0&0&0&0&0&0&0&0\\0&0&1&1&1&1&1&1&2&2\\\end{array}$}\\
\evnrow \hypertarget{desc:327}{}\hyperlink{ps:327}{$327$}&$\MW{4}{5}$&$256$&$16$&$57$&$425$&\cimatrix{$\begin{array}{ccccccc|c}
1&1&1&0&0&0&0&1\\0&0&0&1&1&1&1&2\\\end{array}$}\\
\oddrow \hypertarget{desc:326}{}\hyperlink{ps:326}{$326$}&$\MW{4}{8}$&$320$&$10$&$69$&$525$&\cimatrix{$\begin{array}{cccccccc|cc}
1&1&1&0&0&0&0&0&1&0\\0&0&0&1&1&1&1&1&1&2\\\end{array}$}\\
\evnrow \hypertarget{desc:189}{}\hyperlink{ps:189}{$189$}&$\MW{4}{7}$&$320$&$12$&$69$&$525$&\cimatrix{$\begin{array}{cccccccc|cc}
1&1&1&1&0&0&0&0&1&1\\0&0&0&0&1&1&1&1&1&1\\\end{array}$}\\
\end{longtable}

%% file: table_3.tex
\begin{longtable}{ccccccc}
\caption{Certain $4$-dimensional Fano manifolds with Fano index $r=3$ that arise as complete intersections in toric Fano manifolds.}
\label{tab:index_3}\\
\toprule
\multicolumn{1}{c}{Period ID}&\multicolumn{1}{c}{Description}&\multicolumn{1}{c}{$(-K_X)^4$}&\multicolumn{1}{c}{$\chi(T_X)$}&\multicolumn{1}{c}{$h^0(-K_X)$}&\multicolumn{1}{c}{$h^0(-2K_X)$}&\multicolumn{1}{c}{Weights and line bundles}\\
\midrule
\endfirsthead
\multicolumn{7}{l}{\tiny Continued from previous page.}\\
\addlinespace[1.7ex]
\midrule
\multicolumn{1}{c}{Period ID}&\multicolumn{1}{c}{Description}&\multicolumn{1}{c}{$(-K_X)^4$}&\multicolumn{1}{c}{$\chi(T_X)$}&\multicolumn{1}{c}{$h^0(-K_X)$}&\multicolumn{1}{c}{$h^0(-2K_X)$}&\multicolumn{1}{c}{Weights and line bundles}\\
\midrule
\endhead
\midrule
\multicolumn{7}{r}{\tiny Continued on next page.}\\
\endfoot
\bottomrule
\endlastfoot
\oddrow \hypertarget{desc:37}{}\hyperlink{ps:37}{$37$}&$\FI{4}{3}$&$243$&$27$&$55$&$406$&\cimatrix{$\begin{array}{cccccc|c}
1&1&1&1&1&1&3\\\end{array}$}\\
\evnrow \hypertarget{desc:33}{}\hyperlink{ps:33}{$33$}&$\FI{4}{4}$&$324$&$12$&$70$&$532$&\cimatrix{$\begin{array}{ccccccc|cc}
1&1&1&1&1&1&1&2&2\\\end{array}$}\\
\oddrow \hypertarget{desc:20}{}\hyperlink{ps:20}{$20$}&$\PP^2\times \PP^2$&$486$&$9$&$100$&$784$&\cimatrix{$\begin{array}{cccccc}
1&1&1&0&0&0\\0&0&0&1&1&1\\\end{array}$}\\
\end{longtable}

%% file: table_4.tex
\begin{longtable}{ccccccc}
\caption{The $4$-dimensional Fano manifold with Fano index $r=4$ as a hypersurface in a toric Fano manifold.}
\label{tab:index_4}\\
\toprule
\multicolumn{1}{c}{Period ID}&\multicolumn{1}{c}{Description}&\multicolumn{1}{c}{$(-K_X)^4$}&\multicolumn{1}{c}{$\chi(T_X)$}&\multicolumn{1}{c}{$h^0(-K_X)$}&\multicolumn{1}{c}{$h^0(-2K_X)$}&\multicolumn{1}{c}{Weights and line bundles}\\
\midrule
\endfirsthead
\multicolumn{7}{l}{\tiny Continued from previous page.}\\
\addlinespace[1.7ex]
\midrule
\multicolumn{1}{c}{Period ID}&\multicolumn{1}{c}{Description}&\multicolumn{1}{c}{$(-K_X)^4$}&\multicolumn{1}{c}{$\chi(T_X)$}&\multicolumn{1}{c}{$h^0(-K_X)$}&\multicolumn{1}{c}{$h^0(-2K_X)$}&\multicolumn{1}{c}{Weights and line bundles}\\
\midrule
\endhead
\midrule
\multicolumn{7}{r}{\tiny Continued on next page.}\\
\endfoot
\bottomrule
\endlastfoot
\oddrow \hypertarget{desc:3}{}\hyperlink{ps:3}{$3$}&$Q^4$&$512$&$6$&$105$&$825$&\cimatrix{$\begin{array}{cccccc|c}
1&1&1&1&1&1&2\\\end{array}$}\\
\end{longtable}

%% file: table_5.tex
\begin{longtable}{ccccccc}
\caption{The $4$-dimensional Fano manifold with Fano index $r=5$ is toric.}
\label{tab:index_5}\\
\toprule
\multicolumn{1}{c}{Period ID}&\multicolumn{1}{c}{Description}&\multicolumn{1}{c}{$(-K_X)^4$}&\multicolumn{1}{c}{$\chi(T_X)$}&\multicolumn{1}{c}{$h^0(-K_X)$}&\multicolumn{1}{c}{$h^0(-2K_X)$}&\multicolumn{1}{c}{Weights and line bundles}\\
\midrule
\endfirsthead
\multicolumn{7}{l}{\tiny Continued from previous page.}\\
\addlinespace[1.7ex]
\midrule
\multicolumn{1}{c}{Period ID}&\multicolumn{1}{c}{Description}&\multicolumn{1}{c}{$(-K_X)^4$}&\multicolumn{1}{c}{$\chi(T_X)$}&\multicolumn{1}{c}{$h^0(-K_X)$}&\multicolumn{1}{c}{$h^0(-2K_X)$}&\multicolumn{1}{c}{Weights and line bundles}\\
\midrule
\endhead
\midrule
\multicolumn{7}{r}{\tiny Continued on next page.}\\
\endfoot
\bottomrule
\endlastfoot
\oddrow \hypertarget{desc:1}{}\hyperlink{ps:1}{$1$}&$\PP^4$&$625$&$5$&$126$&$1001$&\cimatrix{$\begin{array}{ccccc}
1&1&1&1&1\\\end{array}$}\\
\end{longtable}